\documentclass[a4paper,12pt]{article}

\usepackage{amsmath}
\usepackage{amsthm}
\usepackage{amssymb}
\usepackage[utf8]{inputenc}
\usepackage{graphicx}

\def\Real{{\mathbb R}}
\def\Cmpx{{\mathbb C}}
\def\Intg{{\mathbb Z}}
\def\Natn{{\mathbb N}}

\def\union{\bigcup}

\def\norm#1{\left|#1\right|}
\def\innerprod(#1,#2){{\left<#1\,,\,#2\right>}}
\def\Set#1{{\left\{#1\right\}}}
\def\qquadtext#1{\qquad\textup{#1}\qquad}
\def\qquadand{\qquadtext{and}}
\def\quadtext#1{\quad\textup{#1}\quad}
\def\quadand{\quadtext{and}}

\def\dfrac#1#2{\frac{d #1}{d #2}}
\def\ddfrac#1#2{\frac{d^2 #1}{d {#2}^2}}

\def\FunSet{{\cal F}}
\def\FunSetO{{\cal F}_0}
\def\ColA{\boldsymbol A}
\def\ColB{\boldsymbol B}
\def\ColC{\boldsymbol C}
\def\ColR{\boldsymbol R}
\def\ColN{\boldsymbol N}
\def\ColO{\boldsymbol 0}
\def\AlF{{\cal H}}
\def\Gen{{\cal G}}
\def\ModO{{\cal E}^{(0)}}
\def\Mod{{\cal E}}
\def\Intinf{{\int_{-\infty}^\infty}}

\def\Ord{{\cal O}}
\def\fbar{{\overline f}}
\def\ftilde{{\tilde f}}
\def\AA{{\mathbb{A}}}
\def\NullO{{\cal N}^{(0)}}
\def\Null{{\cal N}}

\begin{document}

\title{Colombeau Algebra: A pedagogical introduction}
\author{Jonathan Gratus\thanks{j.gratus@lancaster.ac.uk}\\Physics Department, Lancaster University,
  Lancaster LA1 4YB\\and the Cockcroft Institute.}

\maketitle

\begin{abstract}
A simple pedagogical introduction to the Colombeau algebra of
generalised functions is presented, leading the standard definition.
\end{abstract}

\section{Introduction}
\label{ch_Intr}

This is a pedagogical introduction to the Colombeau algebra of
generalised functions. I will limit myself to the Colombeau Algebra
over $\Real$. Rather than $\Real^n$. This is mainly for clarity. Once
the general idea has been understood the extension to $\Real^n$ is not
too difficult. In addition I have limited the introduction to $\Real$
valued generalised functions. To replace with $\Cmpx$ valued
generalised functions is also rather trivial.

I hope that this guide is useful in your understanding of Colombeau
Algebras. Please feel free to contact me.

There is much general literature on Colombeau Algebras but I found the
books by Colombeau himself\cite{colombeau2011elementary} and the Masters thesis by T\d{a} Ng\d{o}c Tr\'\i
\cite{Tri2005colombeau} useful.
 
\section{Test functions and Distributions}

The set of infinitely differentiable functions on $\Real$ is given by
\begin{align}
\FunSet(\Real)
= 
\Set{\phi:\Real\to\Real\ |\ \phi^{(n)}\text{ exists for all
    $n$}}
\label{Intr_def_FunSet}
\end{align}
Test functions are those function which in addition to being smooth are
zero outside an interval, i.e.
\begin{align}
\FunSetO(\Real) = \Set{\phi\in\FunSet(\Real)
\ |\ \text{there exists $a,b\in\Real$
    such that $f(x)=0$ for $x<a$ and $x>b$}}
\label{Intr_def_FunSetO}
\end{align}

I will assume the reader is familiar with distributions, either in the
notation of integrals or as linear functionals. Thus the most important
distributions is the Dirac-$\delta$. This is defined either as a
``function'' $\delta(x)$ such that
\begin{align}
\Intinf  \delta(x) \phi(x) d x  = \phi(0)
\label{Intr_def_delta}
\end{align}
Or as a distribution $\Delta:\FunSetO(\Real)\to\Real$,
\begin{align}
\Delta[\phi]=\phi(0)
\label{Intr_def_Delta}
\end{align}
We will refer to (\ref{Intr_def_delta}) as the integral notation and
(\ref{Intr_def_Delta}) as the Schwartz notation.
An arbitrary distribution will be written either as $\psi(x)$ for the
integral notation  or $\Psi$ for the Schwartz notation.

\section{Function valued distributions}
The first step in understanding the Colombeau Algebra is to convert
distributions into a new object which takes a test functions $\phi$
and gives a functions 
\begin{align*}
\ColA:\FunSetO(\Real)\to\FunSet(\Real)
\end{align*}
This is achieved by using translation of the test functions. Given
$\phi\in\FunSetO$ then let 
\begin{align}
\phi^y\in\FunSetO(\Real)
\,,\qquad
\phi^y(x)=\phi(x-y)
\label{Intr_phi_y}
\end{align}
Then in integral notation
\begin{align}
\overline{\psi}[\phi](y) = \Intinf  \psi(x)\phi(x-y) dx
\label{Intr_def_overline_psi}
\end{align}
and in Schwartz notation
\begin{align}
\overline{\Psi}[\phi](y) = \Psi[\phi^y]
\label{Intr_def_overline_Psi}
\end{align}
We will define the Colombeau Algebra in such a way that they include
the elements $\overline{\psi}$ and $\overline{\Psi}$.
The overline will be used to covert distributions into elements of the
Colombeau algebra.

We label the set of all function valued functionals 
\begin{align}
\AlF(\Real)=\Set{\ColA:\FunSetO(\Real)\to\FunSet(\Real)}
\label{Intr_def_AlF}
\end{align}
We see below that we need to restrict $\AlF(\Real)$ further in order
to define the Colombeau algebra $\Gen(\Real)$.

Observe that we use a slightly non standard notation.  Here
$\ColA[\phi]:\Real\to\Real$ is a function, so that given a point
$x\in\Real$ then $\ColA[\phi](x)\in\Real$. One can equally write
$\ColA[\phi](x)=\ColA(\phi,x)$, which is the standard notation in the
literature. However I claim that the notation
$\ColA[\phi](x)$ does have advantages.

\section{Three special examples.} 
\begin{figure}[t]
\begin{center}
\setlength{\unitlength}{0.5\textwidth}
\includegraphics[height=1\unitlength]{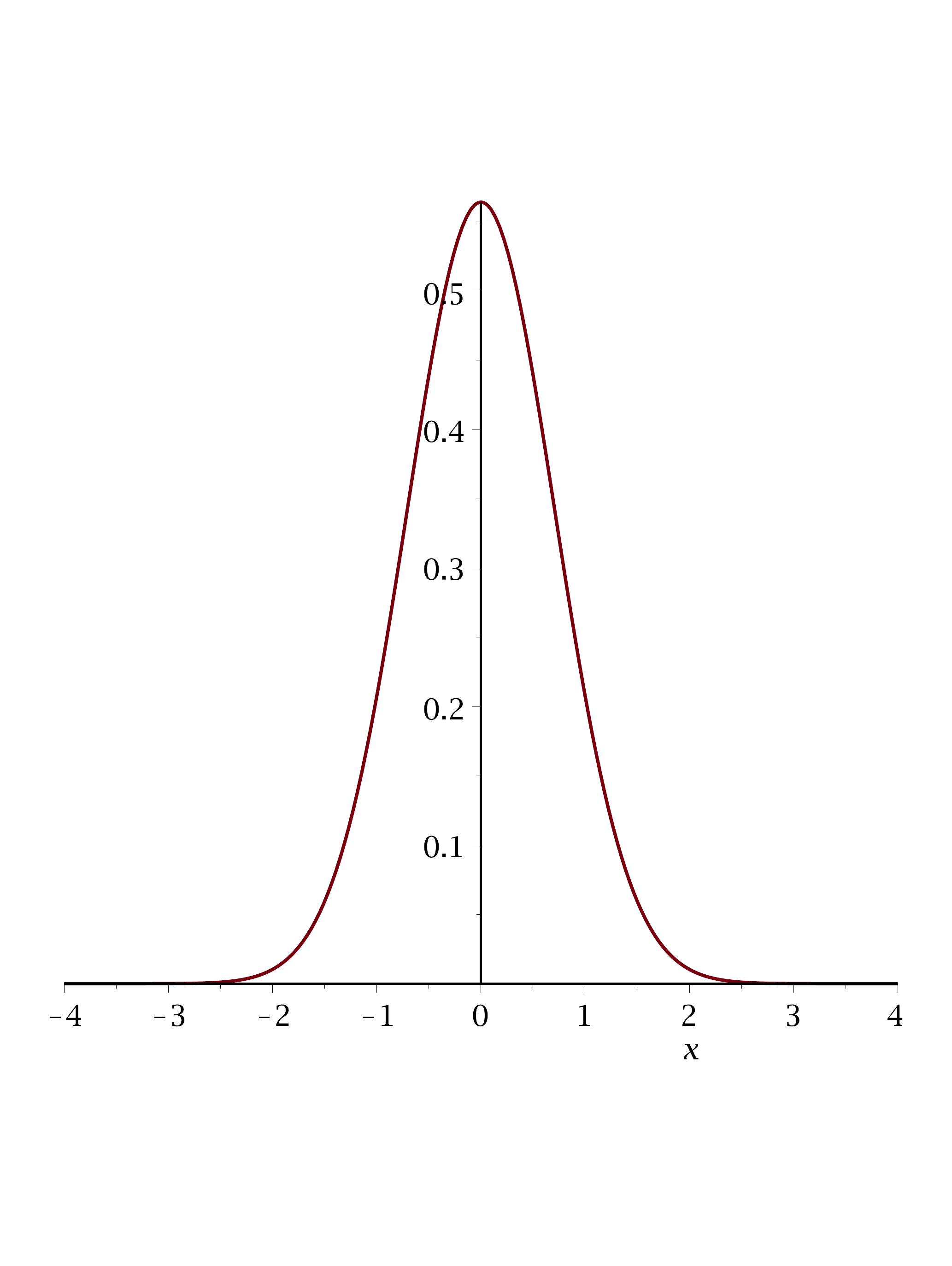}
\qquad\qquad
\includegraphics[height=1\unitlength]{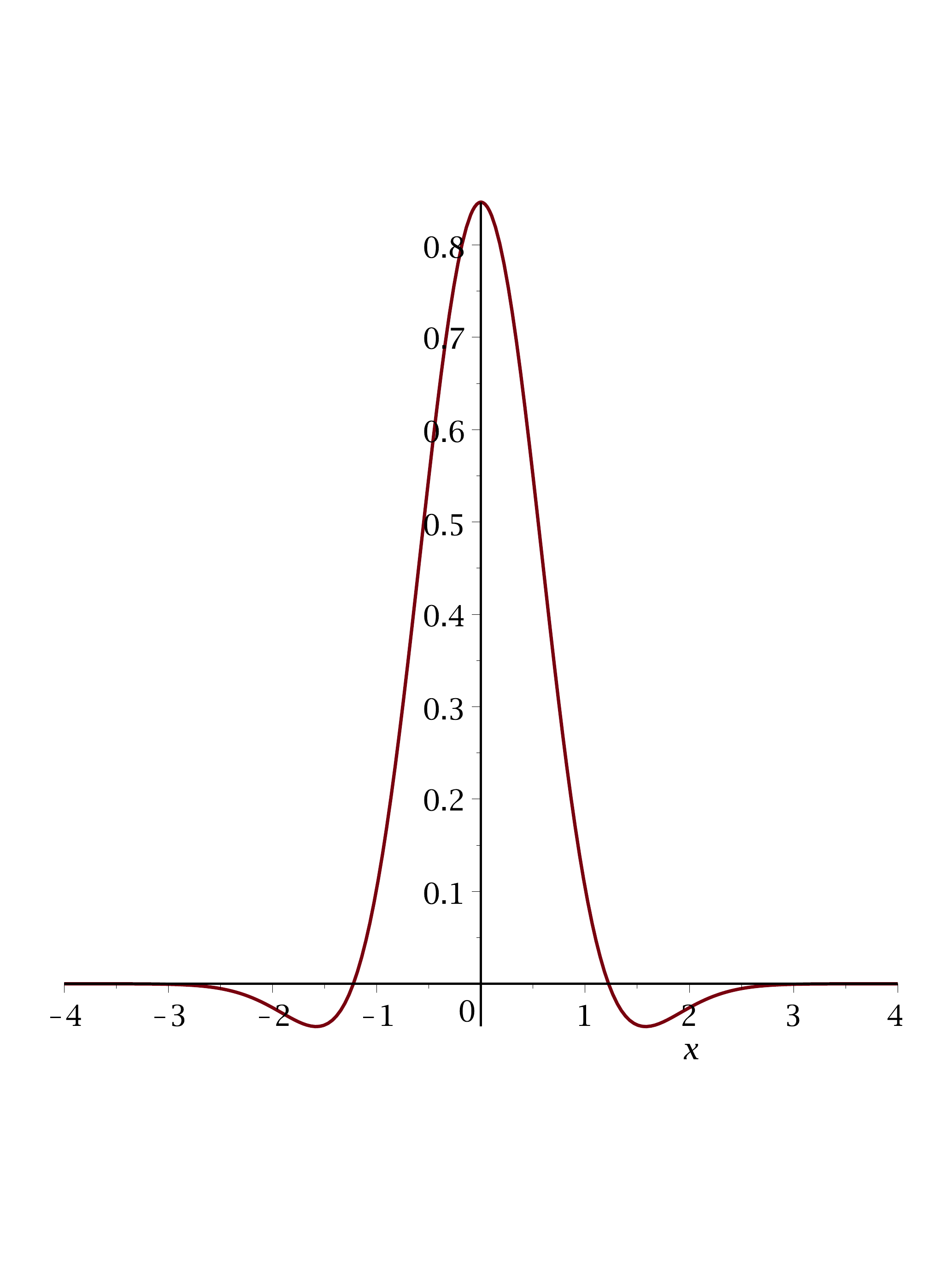}
\end{center}
\vspace{-4em}
\caption{Plots of $\phi_1\in\AA_1$ and $\phi_3\in\AA_3$}
\label{fig_phi}
\end{figure}

For the Dirac-$\delta$ we see that 
\begin{align}
\overline\delta=\overline\Delta=\ColR
\label{Intr_overline_delta}
\end{align}
where $\ColR\in\AlF(\Real)$ is the reflection map
\begin{align}
\ColR[\phi](y)=\phi(-y)
\label{Intr_def_R}
\end{align}
This is because 
\begin{align*}
\overline\delta[\phi](y)
=
\Intinf  \delta(x)\phi(x-y) d x
=
\phi(-y)
\end{align*}
and is Schwartz notation
\begin{align*}
\overline\Delta[\phi](y)
=
\Delta[\phi^y]
=
\phi^y(0)
=
\phi(-y)
\end{align*}

Regular distribution: Given any function $f\in\FunSet$ then there is a
distribution $f^D$ given by
\begin{align}
f^D[\phi]=\Intinf  f(x)\phi(x) dx
\label{Intr_def_fD}
\end{align}
Thus we set $\fbar=\overline{f^D}\in\AlF(\Real)$ as
\begin{align}
\fbar[\phi](y)
&=
{f^D}[\phi^y]
=
\Intinf  f(x) \phi(x-y) dx
\label{Intr_def_overline_f}
\end{align}

The other important generalised functions are the regular
functions. That is given $f\in\FunSet$ we set
\begin{align}
\tilde{f}\in\AlF(\Real)
\,,\qquad
\tilde{f}[\phi]=f
\quadtext{that is}
\tilde{f}[\phi](y)=f(y)
\label{Intr_tilde_f}
\end{align}

\begin{figure}[t]
\begin{center}
\setlength{\unitlength}{0.5\textwidth}
\includegraphics[height=1\unitlength]{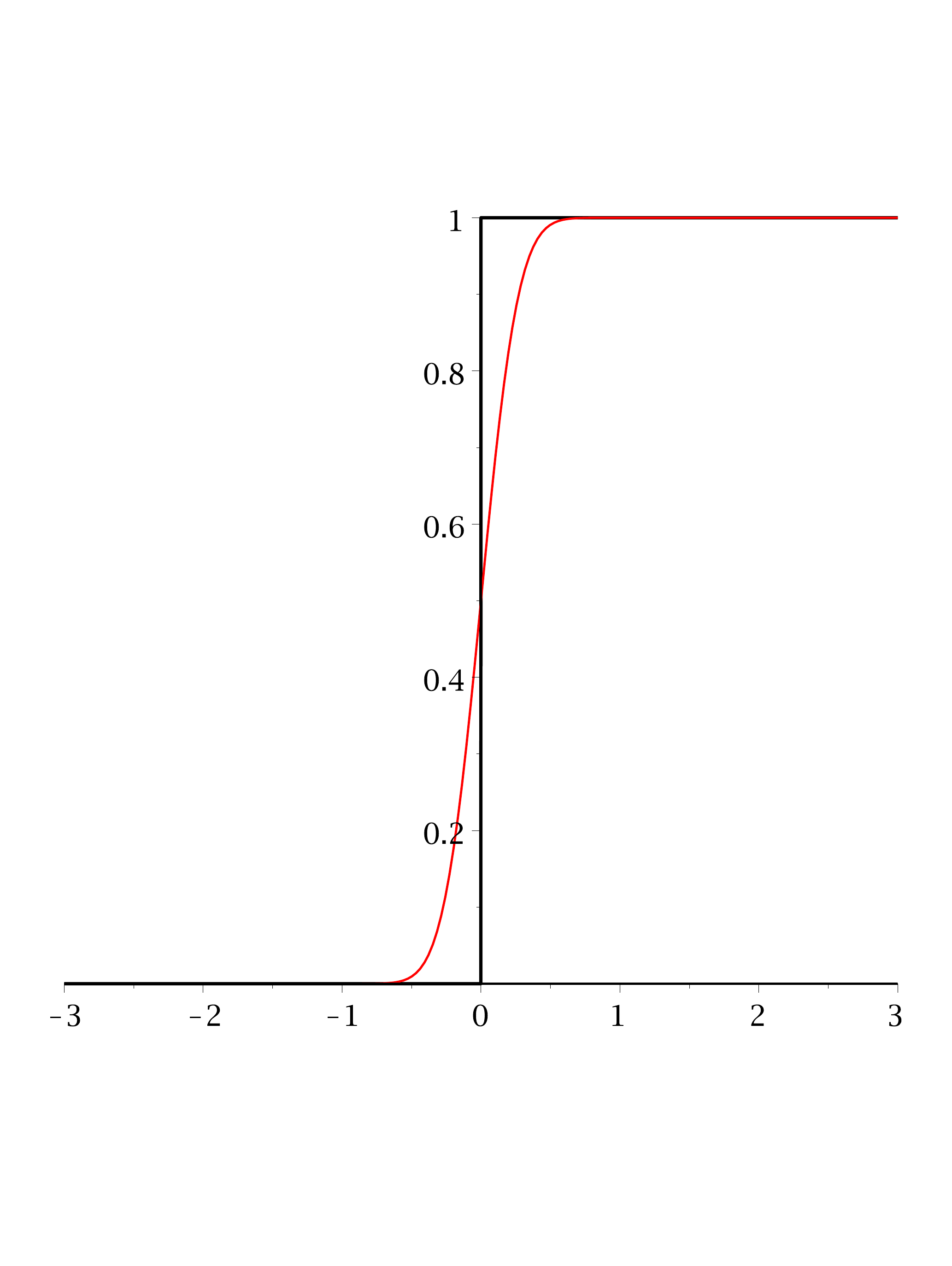}\quad
\includegraphics[height=1\unitlength]{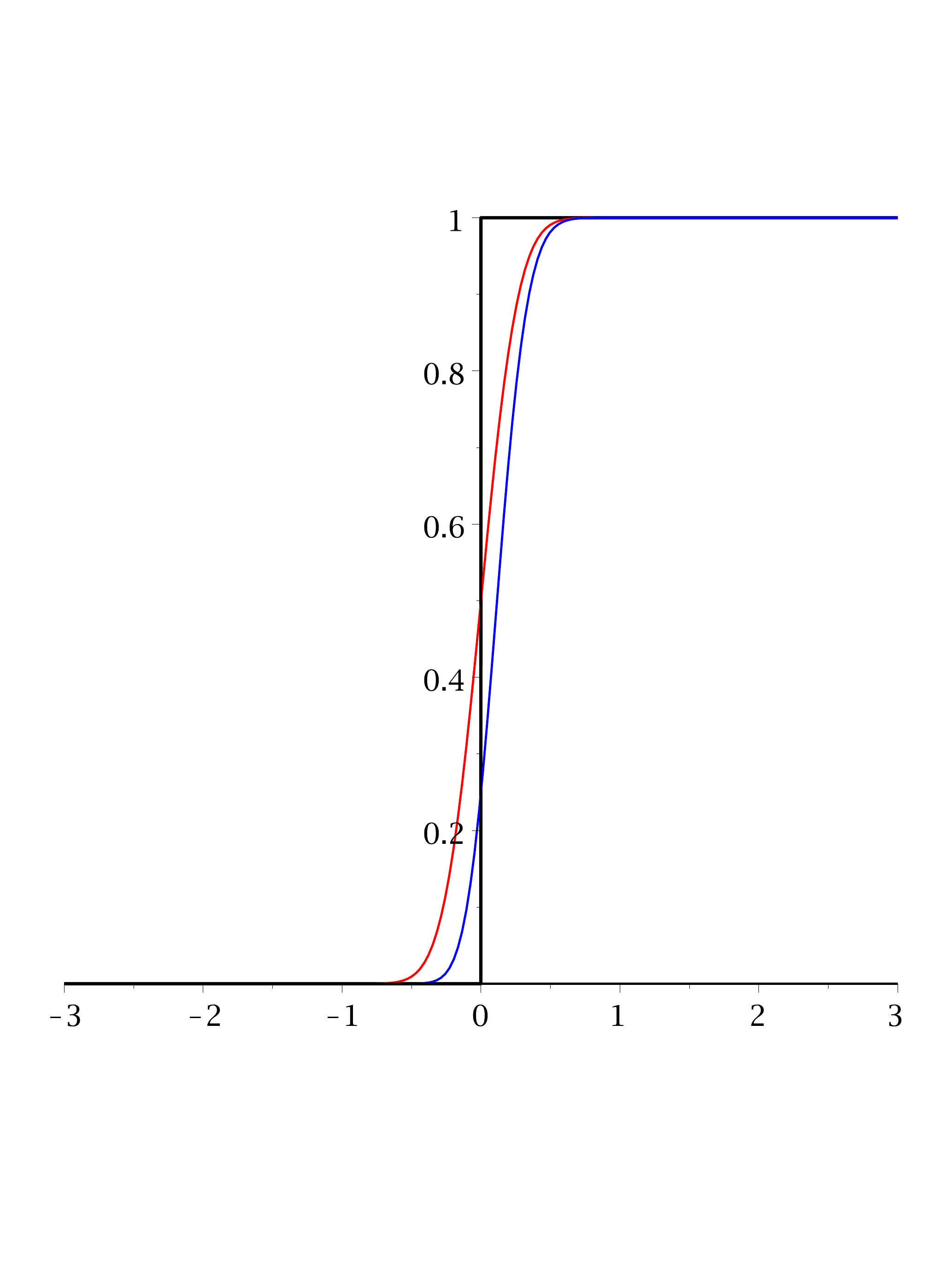}
\end{center}
\vspace{-4em}
\caption{Heaviside (black) and  $\overline\theta[\phi]$ (red) and
  $\big(\overline\theta[\phi]\big)^2$ (blue)}
\label{fig_Heaviside}
\end{figure}

The effect of replacing  $\overline{\psi}[\phi_\epsilon]$ is to smooth out
$\psi$. Examples of $\phi$ are given in figure \ref{fig_phi}. The
action $\overline\theta[\phi]$ where $\theta$ is the Heaviside
function is given in figure \ref{fig_Heaviside}.

\section{Sums and Products}

Given two Generalised functions $\ColA,\ColB\in\AlF(\Real)$ then we
can define there sum and product in the natural way
\begin{align}
\ColA+\ColB\in\AlF(\Real)
\quadtext{via}
(\ColA+\ColB)[\phi]=\ColA[\phi]+\ColB[\phi]
\quadtext{i.e.}
(\ColA+\ColB)[\phi](y)=\ColA[\phi](y)+\ColB[\phi](y)
\label{Intr_def_sum}
\end{align}
and
\begin{align}
\ColA\ColB\in\AlF(\Real)
\quadtext{via}
(\ColA\ColB)[\phi]=\ColA[\phi]\ColB[\phi]
\quadtext{i.e.}
(\ColA\ColB)[\phi](y)=\ColA[\phi](y)\ColB[\phi](y)
\label{Intr_def_prod}
\end{align}

We see that the product of delta functions
$\overline{\delta}{}^2\in\AlF(\Real)$ is clearly defined. That is
\begin{align*}
\overline{\delta}{}^2[\phi](y)
=
(\overline{\delta}[\phi]\overline{\delta}[\phi])(y)
=
\overline{\delta}[\phi](y)\overline{\delta}[\phi](y)
=
\big(\phi(-y)\big)^2
\end{align*}
Although this is a generalised function, it does not correspond to a
distribution, via (\ref{Intr_def_overline_Psi}). That is there is no
distribution $\Psi$ such that $\overline{\Psi}=(\overline{\delta})^2$.

Likewise we can see from figure \ref{fig_Heaviside} that
$(\overline\theta)^2[\phi]=\big(\overline\theta[\phi]\big)^2
\ne\overline\theta[\phi]$.

\section{Making $\fbar$ and $\ftilde$ equivalent}

Now compare the generalised function $\fbar$ and $\ftilde$
(\ref{Intr_def_overline_f}),(\ref{Intr_tilde_f}). We would like these
two generalised functions to be equivalent, so that we can write
$\fbar\sim\ftilde$. One of the results of making $\fbar\sim\ftilde$ is
that if $f,g\in\FunSet$ then
\begin{align*}
\overline{(fg)}\sim \widetilde{(fg)} = \widetilde{f}\ \widetilde{g} \sim
\overline{f}\ \overline{g}
\end{align*}
In the Colombeau algebra, which is a quotient of
equivalent generalised functions, we say that $\fbar$ and $\ftilde$
are the same generalised function. 

The goal therefore is to restrict the set
of possible $\phi$ so that when they are acted upon by
$(\fbar-\ftilde)$ the difference is {\em small}, where {\em small}
will be made technically precise.
When we think of quantities being small, we need a 1-parameter family
of such quantities such that in the limit the difference
vanishes. Here we label the parameter $\epsilon$ and we are interested
in the limit $\epsilon\to 0$ from above, i.e. with $\epsilon>0$.
Given a one parameter set of functions $g_\epsilon\in\FunSet$ then one
meaning to say
$g_\epsilon$ is small is if $g_\epsilon(y)\to 0$ for all $y$. However
we would like a whole hierarchy of smallness. That is for any
$q\in\Natn_0=\Natn\union\Set{0}$ then we can say
\begin{align}
g_\epsilon=\Ord(\epsilon^q)
\label{Intr_g_def_Ord}
\end{align}
if $\epsilon^{-q}g_\epsilon(y)$ is bounded as $\epsilon\to 0$. Note
that we use bounded, rather that tends to zero. However, clearly, 
if $g_\epsilon=\Ord(\epsilon^q)$ then $\epsilon^{-q+1}g_\epsilon\to 0$
as $\epsilon\to 0$. 

We will also need the notion of
$g_\epsilon=\Ord(\epsilon^q)$ where $q<0$. Thus we wish to consider
functions which blow up as $\epsilon\to 0$, but not too quickly. Such 
functions will be called {\em moderate}.

Technically we say $g_\epsilon$ satisfies (\ref{Intr_g_def_Ord}) if
for any interval $(a,b)$ there exists $C>0$ and $\eta>0$ such that 
\begin{align}
\epsilon^{-q}|g_\epsilon(x)|<C
\qquadtext{for all}
a\le y\le b
\qquadand
0<\epsilon<\eta
\label{Intr_def_Ord}
\end{align}

We introduce the parameter $\epsilon$ via the test functions,
replacing $\phi\in\FunSetO$ with $\phi_\epsilon\in\FunSetO$ where
\begin{align}
\phi_\epsilon(x)=\frac{1}{\epsilon} \phi\Big(\frac{x}{\epsilon}\Big)
\label{Intr_def_phi_epsilon}
\end{align}
Observe at as $\epsilon\to 0$ then $\phi_\epsilon$ becomes narrower
and taller, in a definite sense more like a $\delta$-function.  Thus we
consider a generalised function $\ColA$ to be small, if for some
appropriate set of test functions $\phi\in\FunSetO$ and for some
$q\in\Intg$, $\ColA[\phi_\epsilon]=\Ord(\epsilon^q)$.

\begin{figure}
\setlength{\unitlength}{0.24\textwidth}
\begin{tabular}{cc}
\includegraphics[width=2\unitlength]{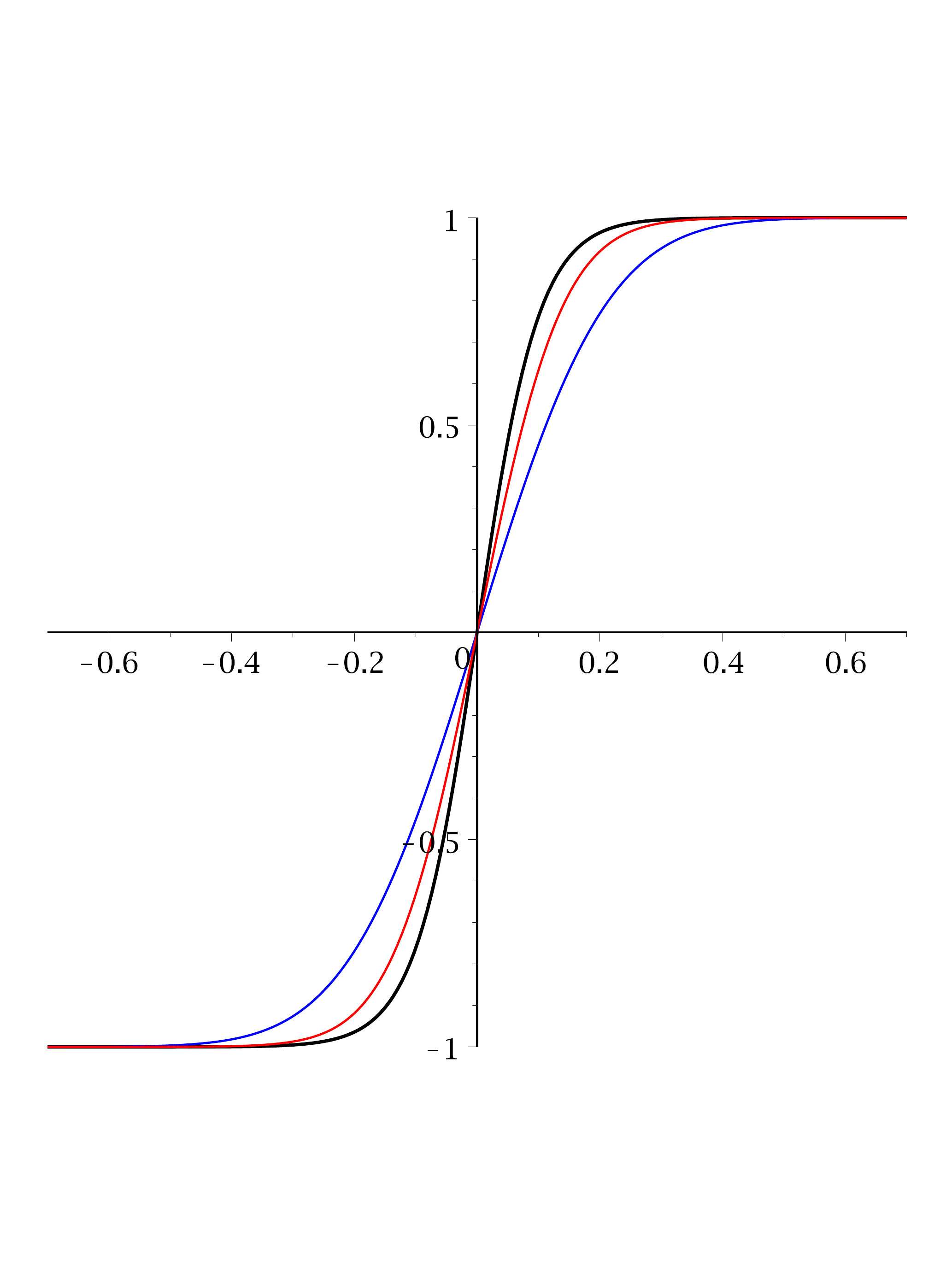}
&
\includegraphics[width=2\unitlength]{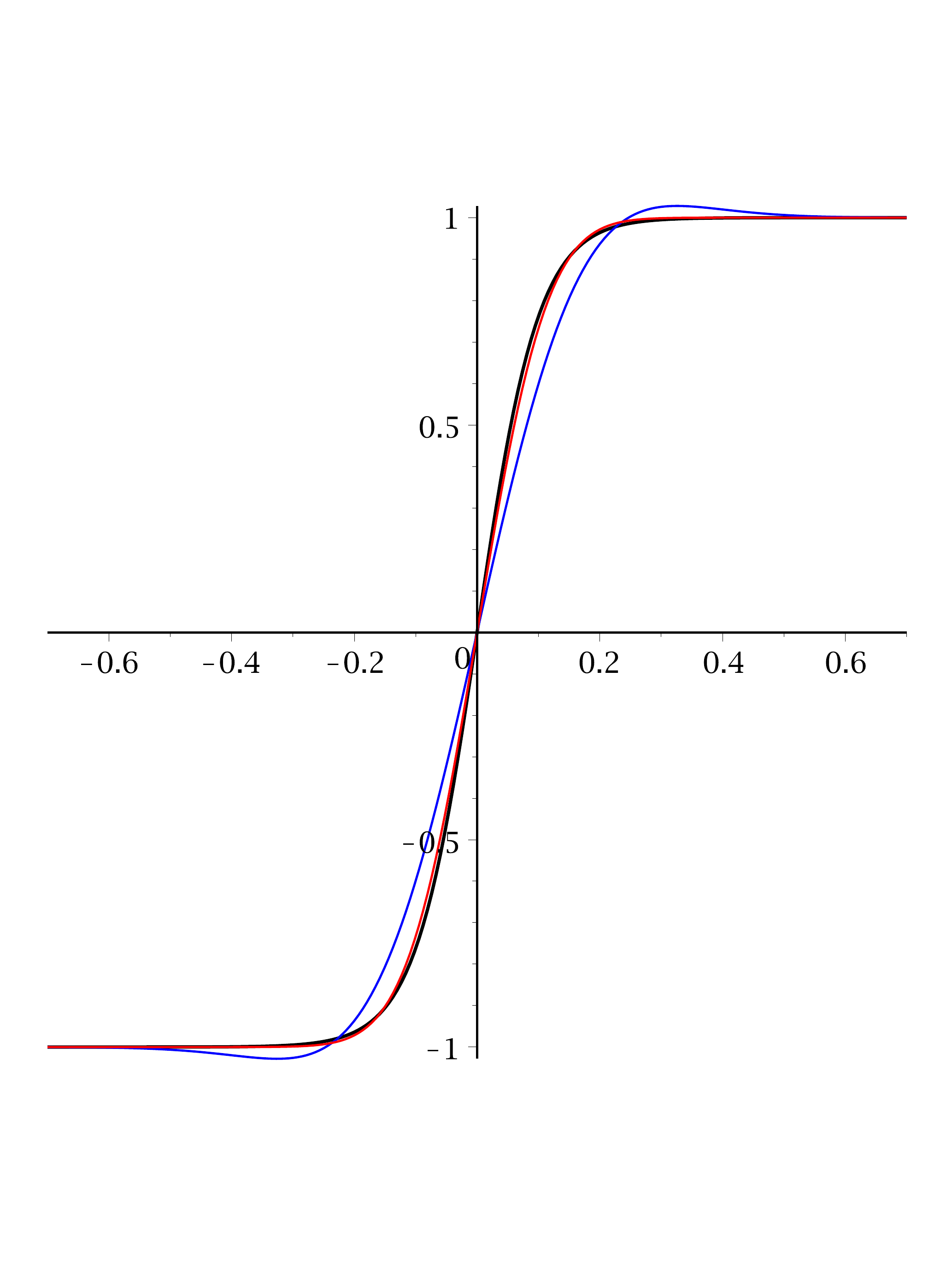}
\\[-2em]
\parbox{2\unitlength}{ $f$ (back),
  $\fbar[\phi_1|_{\epsilon=0.2}]$ (blue) and
  $\fbar[\phi_1|_{\epsilon=0.1}]$ (red).}
&
\parbox{2\unitlength}{ $f$ (back),
  $\fbar[\phi_3|_{\epsilon=0.2}]$ (blue) and
  $\fbar[\phi_3|_{\epsilon=0.1}]$ (red).}
\\[-2em]
\includegraphics[width=2\unitlength]{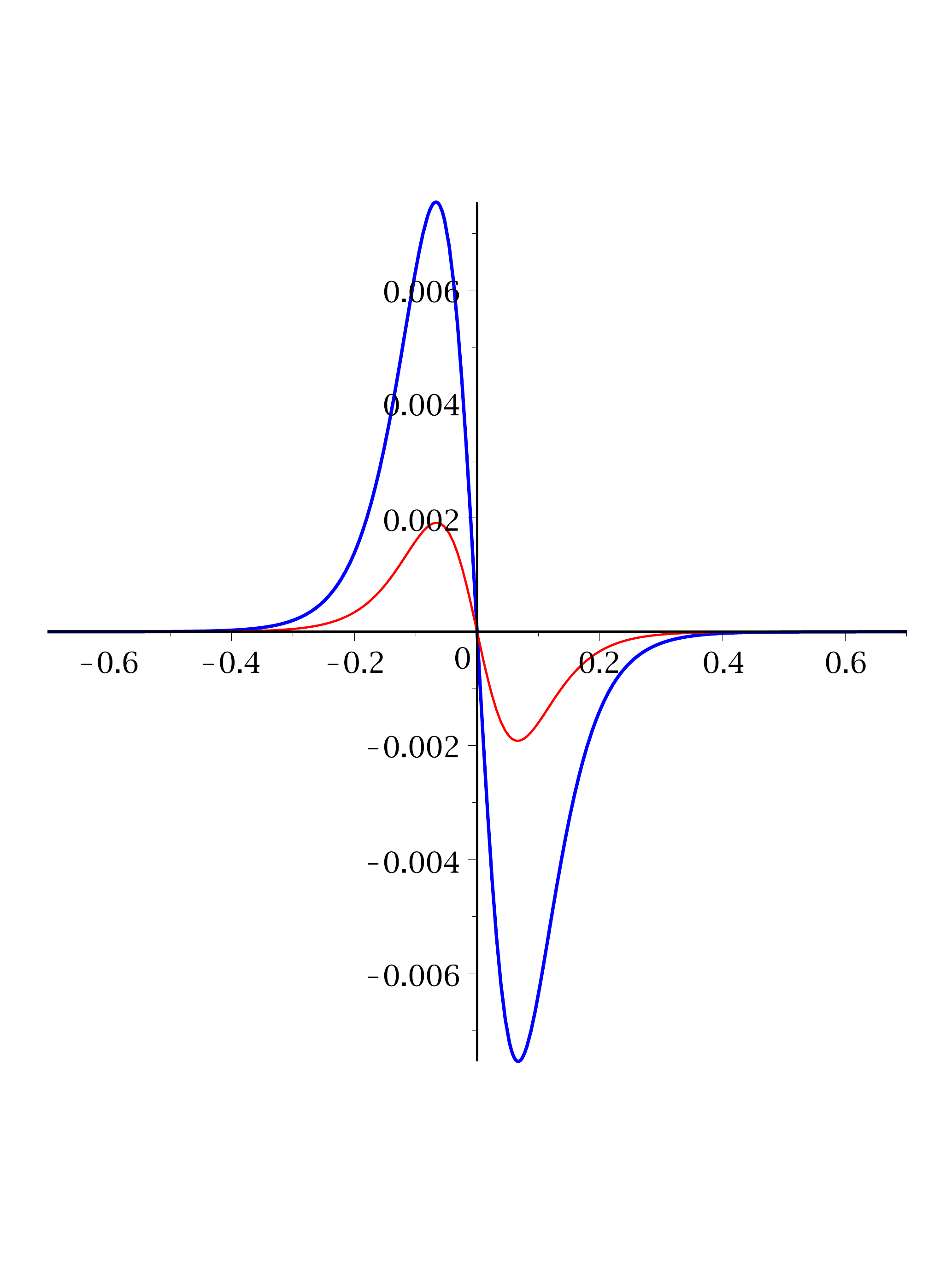}
&
\includegraphics[width=2\unitlength]{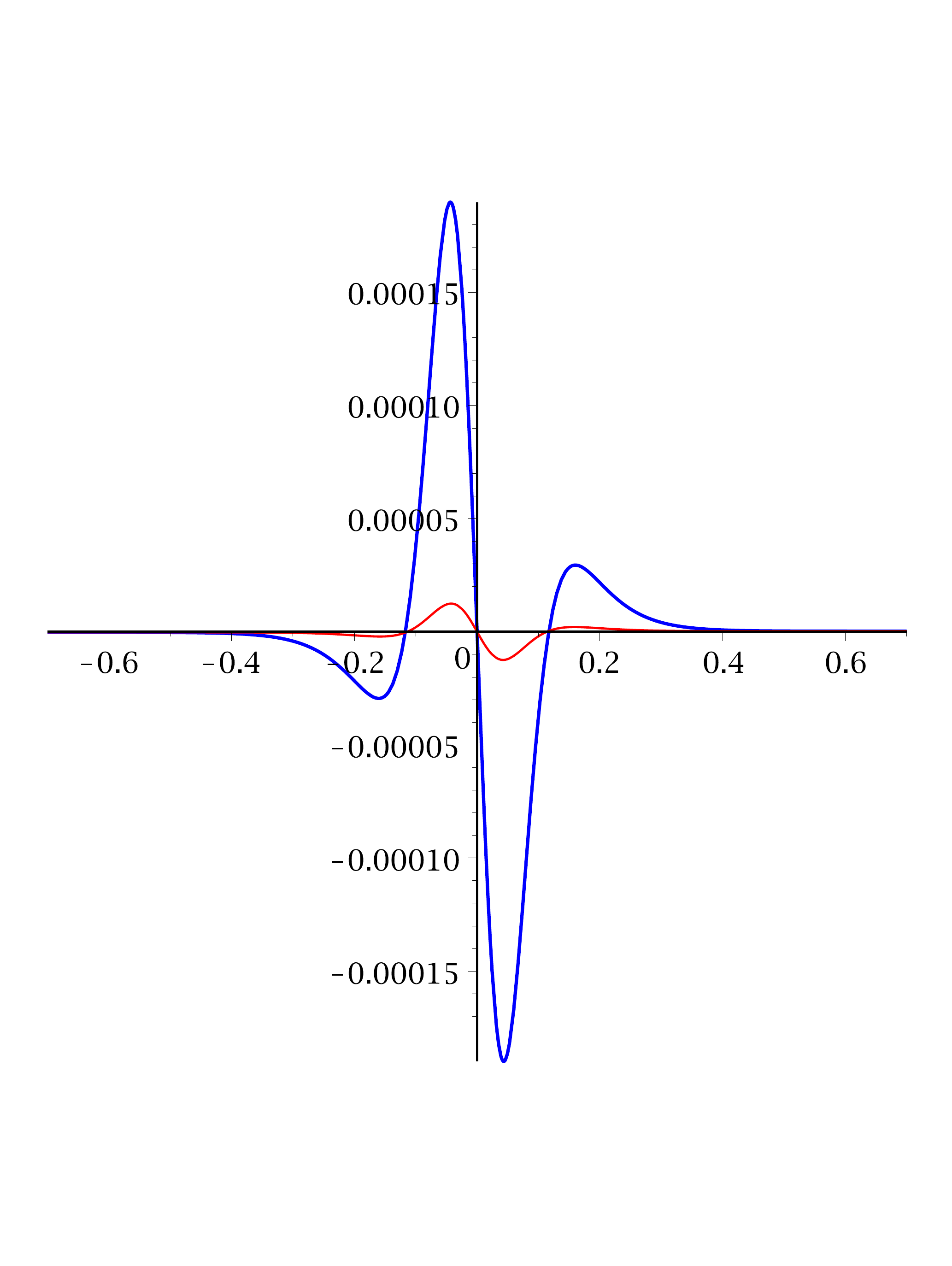}
\\[-2em]
\parbox{2\unitlength}{
  $(\fbar-\ftilde)[\phi_1|_{\epsilon=0.02}]$ (blue) and\\
  $(\fbar-\ftilde)[\phi_1|_{\epsilon=0.01}]$ (red).}
&
\parbox{2\unitlength}{
  $(\fbar-\ftilde)[\phi_3|_{\epsilon=0.02}]$ (blue) and\\
  $(\fbar-\ftilde)[\phi_3|_{\epsilon=0.01}]$ (red).  }
\end{tabular}
\caption{Plots of $\fbar[\phi_{\epsilon}]$ with $f(x)=\tanh(10x)$}
\label{fig_fbar_phi}
\end{figure}

Let us first restrict $\phi\in\FunSetO$ to be test function which
integrate to $1$. That is we define $\AA_0\subset\FunSetO$,
\begin{align}
\AA_0=\Set{\phi\in\FunSetO\,\Big|\,
\Intinf  \phi(x) dx = 1}
\label{Intr_int_phi_1}
\end{align}
Given $\phi\in\AA_0$ and setting $z=(x-y)/\epsilon$ so that $x=y+\epsilon z$
\begin{equation}
\begin{aligned}
\fbar[\phi_\epsilon](y)
&=
{f^D}[\phi^y_\epsilon]
=
\Intinf  f(x) \phi_\epsilon(x-y) dx
=
\frac{1}{\epsilon}\Intinf  f(x)
\phi\Big(\frac{x-y}{\epsilon}\Big) 
dx
\\&=
\Intinf  f(y+\epsilon z)
\phi(z) 
dz
\end{aligned}
\label{Intr_fbar_int}
\end{equation}
Thus as $\epsilon\to 0$ then $f(y+\epsilon z)\approx f(y)$ so that,
since $\phi\in\AA_0$,
\begin{align*}
\fbar[\phi_\epsilon](y)
=
\Intinf  f(y+\epsilon z)
\phi(z) 
dz
\approx
\Intinf  f(y)
\phi(z) 
dz
=
f(y)
\Intinf  
\phi(z) 
dz
=
f(y)
=
\ftilde[\phi_\epsilon](y)
\end{align*}
In fact since $f(y+\epsilon z)-f(y)=\Ord(\epsilon)$ we can show using
(\ref{Intr_def_Ord}) that
\begin{align}
\text{if}\qquad
\phi\in\AA_0
\qquadtext{then}
(\fbar-\ftilde)[\phi_\epsilon]=\Ord(\epsilon)
\label{Intr_fbar_ftilde_epsilon}
\end{align}

This is good so far, but we want to further restrict the set $\phi$ so
that we can satisfy 
\begin{align}
(\fbar-\ftilde)[\phi_\epsilon]=\Ord(\epsilon^q)
\label{Intr_f_diff_Ord(n)}
\end{align}
to any order of $\Ord(\epsilon^q)$. 

Taylor expanding $f(y+\epsilon z)$ to order $\Ord(\epsilon^{q+1})$ we have
\begin{align*}
f(y+\epsilon z)
=
\sum_{r=0}^q \frac{\epsilon^r z^r f^{(r)}(y)}{r!} 
+
\Ord(\epsilon^{q+1})
\end{align*}
Thus 
\begin{equation}
\begin{aligned}
(\fbar-\ftilde)[\phi_\epsilon](y)
&=
\Intinf 
\big(f(y+\epsilon z)-f(y)\big)\phi(z) dz
=
\Intinf
\Big(
\sum^q_{n=1}
\frac{\epsilon^r z^r f^{(r)}(y)}{r!} 
+
\Ord(\epsilon^{q+1})
\Big)\phi(z) dz
\\&=
\sum^q_{n=1} \frac{\epsilon^r f^{(r)}(y)}{r!}
\Intinf
z^r \phi(z) dz
+
\Ord(\epsilon^{q+1})
\end{aligned}
\label{Intr_fbar_f_O(q)}
\end{equation}
Thus we can satisfy (\ref{Intr_g_def_Ord}) to order
$\Ord(\epsilon^{q+1})$ if the first $q$ moments of $\phi(z)$ vanish:
\begin{align*}
\Intinf z^r \phi(z) dz = 0
\quadtext{for}
1\le r \le q
\end{align*}

We now define all the elements with vanishing moments. 
\begin{align}
\AA_q=\Set{\phi\in\FunSetO(\Real)\,\Big|\,
\Intinf \phi(z) dz=1 \quadand
\Intinf z^r \phi(z) dz = 0
\quadtext{for}
1\le r \le q}
\label{Intr_def_AAq}
\end{align}
So clearly $\AA_{q+1}\subset \AA_{q}$.  
We can show that these functions exist.
Thus from (\ref{Intr_fbar_f_O(q)})
we have 
\begin{align}
\phi\in \AA_q \quadtext{implies}
(\fbar-\ftilde)[\phi_\epsilon] = \Ord(\epsilon^{q+1})
\label{Intr_A_q_fbar_f_O(q)}
\end{align}
Two example test functions $\phi_1\in\AA_1$ and $\phi_3\in\AA_3$ are
given in figure \ref{fig_phi}. The result $\fbar[\phi_\epsilon]$,
(\ref{Intr_def_overline_f}), (\ref{Intr_fbar_int}) is given in fig
\ref{fig_fbar_phi}.

The easiest way to construct $\phi\in\AA_q$ is to choose a test function
$\psi$ and then set
\begin{align*}
\phi(z)=\lambda_0\psi(z)+\lambda_1\psi'(z)+\cdots+\lambda_{q-1}\psi^{(q-1)}(z)
\end{align*}
where $\lambda_0,\ldots,\lambda_{q-1}\in\Real$ are constants
determined by (\ref{Intr_def_AAq}).

\section{Null and moderate generalised functions.}

As we stated we wanted $\fbar$ and $\ftilde$ to be considered
equivalent. From (\ref{Intr_A_q_fbar_f_O(q)}) we have $\phi\in\AA_q$
then $(\fbar-\ftilde)[\phi_\epsilon]=\Ord(\epsilon^{q+1})$. We
generalise this notion. We say that $\ColA,\ColB\in\AlF(\Real)$ are
equivalent, $\ColA\sim\ColB$, if for all $q\in\Natn$ there is a
$p\in\Natn$ such that
\begin{align}
\phi\in\AA_p
\quadtext{implies}
\ColA[\phi_\epsilon]-\ColB[\phi_\epsilon]=\Ord(\epsilon^q)
\label{Intr_Ap_Oq}
\end{align}

We label $\NullO(\Real)\subset\AlF(\Real)$ the set of all elements
which are {\em null}, that is equivalent to the zero element
$\ColO\in\AlF(\Real)$ that is
\begin{align*}
\NullO(\Real)
&=
\Set{\ColA\in\AlF(\Real)\,|\, \ColA\sim\ColO}
\end{align*}
I.e.
\begin{align}
\NullO(\Real)
&=
\Set{\ColA\in\AlF(\Real)\,|\,\text{for all } p\in\Natn
\text{ there exists }q\in\Natn
\text{ such that for all }
\phi\in\AA_q,\
\ColA[\phi_\epsilon]=\Ord(\epsilon^p)}
\label{Intr_def_Null0}
\end{align}

Examples of null elements are of course
$\fbar-\ftilde\in\NullO(\Real)$, which is true by
construction. Another example is $\ColN\in\NullO(\Real)$ which is
given by
\begin{align}
\ColN[\phi](y)=\phi(1)
\label{Intr_Col_N}
\end{align}
Since for any $\phi\in\AA_0$ there exists $\eta>0$ such that $1/\eta$
is outside the support of $\phi$. Thus  
$\phi_\epsilon(1)=0$ for all $\epsilon<\eta$ and hence
$\ColN[\phi_\epsilon]=0$ so $\ColN\in\NullO(\Real)$. However, although
$\ColN\in\NullO$, we can choose $\phi$ so that
$\ColN[\phi](y)=\phi(1)$ is any value we choose. Thus knowing that a
generalised function $\ColA$ is null says nothing about the value of
$\ColA[\phi]$ but only the limit of $\ColA[\phi_\epsilon]$ as
$\epsilon\to0$.

We would like $\NullO(\Real)$ to form an ideal in $\AlF(\Real)$, that
is that if $\ColA,\ColB\in\NullO(\Real)$ and $\ColC\in\AlF(\Real)$ then
\begin{itemize}
\item
$\ColA+\ColB\in\NullO(\Real)$ and 
\item
$\ColA\ColC\in\NullO(\Real)$.
\end{itemize}
It is easy to see that the first of these is automatically
satisfied. However the second requires one additional requirement. We
need
\begin{align}
\ColC[\phi_\epsilon]=\Ord(\epsilon^{-N})
\label{Intr_C_O(-N)}
\end{align}
for some $N\in\Intg$. Thus although $\ColC[\phi_\epsilon]\to\infty$ as
$\epsilon\to 0$ we don't want it to blow up to quickly. 
Now we have the following:

Given $\ColA\in\NullO(\Real)$ and $\ColC$ satisfying
(\ref{Intr_C_O(-N)}) and given $q\in\Natn_0$ then there exists
$p\in\Intg$ such that $\phi\in\AA_p$ implies
$\ColA[\phi_\epsilon]=\Ord(\epsilon^{q+N})$.  
Hence
\begin{align*}
(\ColA\ColC)[\phi_\epsilon]
=
\ColA[\phi_\epsilon] \ColC[\phi_\epsilon]
=
\Ord(\epsilon^{q+N})\Ord(\epsilon^{-N})
=
\Ord(\epsilon^{q})
\end{align*}
hence $\ColA\ColC\in\NullO(\Real)$.
We call the set of elements $\ColC\in\AlF(\Real)$ satisfying
(\ref{Intr_C_O(-N)}), {\em moderate} and set of moderate functions
\begin{align}
\ModO(\Real)
=
\Set{\ColA\in\AlF(\Real)\,\big|\,
\text{there exists } N\in\Natn
\text{ such that for all }
\phi\in\AA_0,\
\ColA[\phi_\epsilon]=\Ord(\epsilon^{-N})
}
\label{Intr_def_Mod0}
\end{align}

Examples of moderate functions include
\begin{align*}
\overline\Delta[\phi_\epsilon](y)
=
\phi_\epsilon(-y)
=
\frac{1}{\epsilon}\phi\Big(-\frac{y}{\epsilon}\Big)
=
\Ord(\epsilon^{-1})\,,\qquad
\big(\overline\Delta\big)^n[\phi_\epsilon]=\Ord(\epsilon^{-n})
\end{align*}
and
\begin{align*}
\ftilde[\phi_\epsilon](y)
=
f(y)
=
\Ord(\epsilon^{0})
\end{align*}

\section{Derivatives}

The last part in the construction of the Colombeau Algebra is to
extend all the definitions so that they also apply to the derivatives
$\displaystyle\dfrac{\ColA[\phi]}{y}$,
$\displaystyle\ddfrac{\ColA[\phi]}{y}$, etc. We require that not only
does a moderate function not blow up too quickly, but neither do its
derivatives, i.e.
\begin{align}
(\ColA[\phi])^{(n)}
=
\frac{d^n}{dy^n} \big(\ColA[\phi]\big)
\in\ModO(\Real)
\label{Intr_def_D^n}
\end{align}
Thus we define the set of moderate function as
\begin{align}
\Mod(\Real)
=
\Set{\ColA\in\ModO(\Real)\,\Big|\,
\big(\ColA[\phi]\big)^{(n)}
\in\ModO(\Real)
\text{ for all }
n\in\Natn,\phi\in\AA_0}
\label{Intr_def_ModA}
\end{align}
That is
\begin{equation}
\begin{aligned}
\Mod(\Real)
=
\Big\{\ColA\in\AlF(\Real)\,\Big|\,
\text{for all }
n\in\Natn_0
\text{ there exists }
N\in\Natn
\text{ such that for all }
\phi\in\AA_0,
\qquad&
\\
\big(\ColA[\phi_\epsilon]\big)^{(n)}=\Ord(\epsilon^{-N})
\Big\}&
\end{aligned}
\label{Intr_def_Mod}
\end{equation}
Likewise we require that for two generalised functions to be
equivalent then we require that all the derivatives are small
\begin{align}
\Null(\Real)
=
\Set{\ColA\in\NullO(\Real)\,\Big|\,
\big(\ColA[\phi]\big)^{(n)}
\in\NullO(\Real)
\text{ for all }
n\in\Natn}
\label{Intr_def_NullA}
\end{align}
That is
\begin{equation}
\begin{aligned}
\Null(\Real)
=
\Big\{\ColA\in\AlF(\Real)\,\Big|\,
\text{for all }
n\in\Natn_0
\text{ and }
q\in\Natn
\text{ there exists }
p\in\Natn
\text{ such that}\qquad\qquad&
\\\qquad\qquad\text{for all }
\phi\in\AA_p,
\ 
\big(\ColA[\phi_\epsilon]\big)^{(n)}=\Ord(\epsilon^{q})
\Big\}&
\end{aligned}
\label{Intr_def_Null}
\end{equation}

\section{Quotient Algebra}

We write the Colombeau Algebra as a quotient algebra,
\begin{align}
\Gen(\Real)={\Mod(\Real)}/{\Null(\Real)}
\label{Intr_def_Gen}
\end{align}
This means that, with regard to elements $\ColA,\ColB\in\Mod(\Real)$
we say $\ColA\sim\ColB$ if $\ColA-\ColB\in\Null(\Real)$. For elements
in $\ColA,\ColB\in\Gen(\Real)$ we simply write $\ColA=\ColB$. 

Given $\ColA\in\Gen(\Real)$, then in order to get an actual number we
must first choose a representative $\ColB\in\Mod(\Real)$ of
$\ColA\in\Gen(\Real)$, then we must choose $\phi\in\AA_0$ and
$y\in\Real$ then the quantity $\ColB[\phi](y)\in\Real$.

\section{Summary}

We can summarise the steps needed to go from distributions to
Colombeau functions:
\begin{itemize}
\item
Convert distributions which give a number $\Psi[\phi]$ 
as an answer to functionals $\ColA[\phi]$ which give a function as an
answer.
\item
Construct the sets of test functions $\AA_q$, so that
$\fbar\sim\ftilde$, i.e. $\fbar-\ftilde\in\NullO(\Real)$
\item
Limit the generalised functions to elements of $\ModO(\Real)$ so that
the set $\NullO(\Real)\subset\ModO(\Real)$ is an ideal.
\item
Extend the definitions of $\ModO(\Real)$ and $\NullO(\Real)$ to
$\Mod(\Real)$ and $\Null(\Real)$ so that
they also apply to derivatives.
\item
Define the Colombeau Algebra as the quotient
$\Gen(\Real)=\Mod(\Real)/\Null(\Real)$.

\end{itemize}

The formal definition, we define $\Mod(\Real)$ via
(\ref{Intr_def_Mod}), then $\Null(\Real)$ via (\ref{Intr_def_Null})
and (\ref{Intr_def_AAq}). Then define the Colombeau Algebra
$\Gen(\Real)$ as the
quotient (\ref{Intr_def_Gen}).

\section*{Acknowledgement}

The author is grateful for the support provided by STFC (the
Cockcroft Institute ST/G008248/1) and EPSRC (the Alpha-X project
EP/J018171/1.

\bibliographystyle{plain}
\bibliography{Colombeau}

\begin{thebibliography}{1}

\bibitem{colombeau2011elementary}
Jean~Fran{\c{c}}ois Colombeau.
\newblock {\em Elementary introduction to new generalized functions}.
\newblock Elsevier, 2011.

\bibitem{Tri2005colombeau}
{T\d{a} Ng\d{o}c} {Tr{\'\i}} and Tom~H Koornwinder.
\newblock {The Colombeau Theory of Generalized Functions}.
\newblock Master's thesis, KdV Institute, Faculty of Science, University of
  Amsterdam, The Netherlands, 2005.

\end{thebibliography}

\end{document}